\documentclass [12pt]{article}
\usepackage{amsmath,amssymb}

\begin{document}
\title{ Cramer's rule for some quaternion matrix
equations.}
\author {Kyrchei I. I.\footnote{Pidstrygach Institute for Applied Problems
 of Mechanics and Mathematics of NAS of Ukraine,
str. Naukova 3b, Lviv, Ukraine, 79053, kyrchei@lms.lviv.ua}}

\date{}
 \maketitle


\begin{abstract}
Cramer's rules for some left, right and two-sided quaternion
matrix equations are obtained within the framework of the theory
of the column and row determinants.


\end{abstract}

\section{Introduction} Quaternion matrix equations play important roles in both theoretical studies
and numerical computations of quaternion application disciplines
\cite{ad,ja,zh} and have been studied by many experts
\cite{ti,wa1,wa2} in mainly by means of complex representation of
quaternion matrices. In this paper we obtain Cramer's rule for
some left, right and two-sided quaternion matrix equations within
the framework of the theory of
  new  matrix functionals over the quaternion
skew field (the column and row determinants) introduced in
\cite{ky}. In the first point we shortly cite some provisions from
the theory of the column and row determinants which are necessary
for the following. The theory of the column and row determinants
of a quaternion matrix are considered completely  in \cite{ky}.
 In the second point the Cramer rules
for left, right and two-sided quaternion matrix equations are
obtained.

\section{Elements of the theory of the column and row determinants.}
\newtheorem{theorem}{Theorem}[section]
\newtheorem{definition}{Definition}[section]
\newtheorem{proposition}{Proposition}[section]
\newtheorem{lemma}{Lemma}[section]

Let ${\rm M}\left( {n,{\rm {\mathbb{H}}}} \right)$ be the ring of
$n\times n$ quaternion matrices. By ${\rm {\mathbb{H}}}^{m\times
n}$ denote the set of all $m\times n$ matrices over the quaternion
skew field ${\rm {\mathbb{H}}}$ and by ${\rm
{\mathbb{H}}}_{r}^{m\times n} $ denote its subset of matrices of
rank $r$. Suppose $S_{n}$ is the symmetric group on the set
$I_{n}=\{1,\ldots,n\}$.
\begin{definition}
 The $i$th row determinant of ${\rm {\bf A}}=(a_{ij}) \in {\rm
M}\left( {n,{\mathbb{H}}} \right)$ is defined by
 \[{\rm{rdet}}_{ i} {\rm {\bf A}} =
{\sum\limits_{\sigma \in S_{n}} {\left( { - 1} \right)^{n -
r}{a_{i{\kern 1pt} i_{k_{1}}} } {a_{i_{k_{1}}   i_{k_{1} + 1}}}
\ldots } } {a_{i_{k_{1} + l_{1}}
 i}}  \ldots  {a_{i_{k_{r}}  i_{k_{r} + 1}}}
\ldots  {a_{i_{k_{r} + l_{r}}  i_{k_{r}} }}\] \noindent for all $i
= 1,\ldots,n $. The elements of the permutation $\sigma$ are
indices of each monomial. The left-ordered cycle notation of the
permutation $\sigma$ is written as follows,
\[\sigma = \left(
{i\,i_{k_{1}}  i_{k_{1} + 1} \ldots i_{k_{1} + l_{1}} }
\right)\left( {i_{k_{2}}  i_{k_{2} + 1} \ldots i_{k_{2} + l_{2}} }
\right)\ldots \left( {i_{k_{r}}  i_{k_{r} + 1} \ldots i_{k_{r} +
l_{r}} } \right).\] \noindent The index $i$ opens the first cycle
from the left  and other cycles satisfy the following conditions,
$i_{k_{2}} < i_{k_{3}}  < \ldots < i_{k_{r}}$ and $i_{k_{t}}  <
i_{k_{t} + s} $ for all $t = 2,\ldots,r $ and $s =1,\ldots,l_{t}
$.
\end{definition}

\begin{definition}
The $j$th column determinant
 of ${\rm {\bf
A}}=(a_{ij}) \in {\rm M}\left( {n,{\mathbb{H}}} \right)$ is
defined by
 \[{\rm{cdet}} _{{j}}\, {\rm {\bf A}} =
{{\sum\limits_{\tau \in S_{n}} {\left( { - 1} \right)^{n -
r}a_{j_{k_{r}} j_{k_{r} + l_{r}} } \ldots a_{j_{k_{r} + 1}
i_{k_{r}} }  \ldots } }a_{j\, j_{k_{1} + l_{1}} }  \ldots  a_{
j_{k_{1} + 1} j_{k_{1}} }a_{j_{k_{1}} j}}\] \noindent for all $j
=1,\ldots,n $. The right-ordered cycle notation of the permutation
$\tau \in S_{n}$ is written as follows,
 \[\tau =
\left( {j_{k_{r} + l_{r}}  \ldots j_{k_{r} + 1} j_{k_{r}} }
\right)\ldots \left( {j_{k_{2} + l_{2}}  \ldots j_{k_{2} + 1}
j_{k_{2}} } \right){\kern 1pt} \left( {j_{k_{1} + l_{1}}  \ldots
j_{k_{1} + 1} j_{k_{1} } j} \right).\] \noindent The index $j$
opens  the first cycle from the right  and other cycles satisfy
the following conditions, $j_{k_{2}}  < j_{k_{3}}  < \ldots <
j_{k_{r}} $ and $j_{k_{t}}  < j_{k_{t} + s} $ for all $t =
2,\ldots,r $ and $s = 1,\ldots,l_{t}  $.
\end{definition}

Suppose ${\rm {\bf A}}_{}^{i{\kern 1pt} j} $ denotes the submatrix
of ${\rm {\bf A}}$ obtained by deleting both the $i$th row and the
$j$th column. Let ${\rm {\bf a}}_{.j} $ be the $j$th column and
${\rm {\bf a}}_{i.} $ be the $i$th row of ${\rm {\bf A}}$. Suppose
${\rm {\bf A}}_{.j} \left( {{\rm {\bf b}}} \right)$ denotes the
matrix obtained from ${\rm {\bf A}}$ by replacing its $j$th column
with the column ${\rm {\bf b}}$, and ${\rm {\bf A}}_{i.} \left(
{{\rm {\bf b}}} \right)$ denotes the matrix obtained from ${\rm
{\bf A}}$ by replacing its $i$th row with the row ${\rm {\bf b}}$.

We  note some properties of column and row determinants of a
quaternion matrix ${\rm {\bf A}} = \left( {a_{ij}} \right)$, where
$i \in I_{n} $, $j \in J_{n} $ and $I_{n} = J_{n} = {\left\{
{1,\ldots ,n} \right\}}$.
\begin{proposition}  \cite{ky}
If $b \in {\mathbb{H}}$, then
 $ {\rm{rdet}}_{ i} {\rm {\bf A}}_{i.} \left( {b
\cdot {\rm {\bf a}}_{i.}}  \right) = b \cdot {\rm{rdet}}_{ i} {\rm
{\bf A}}$ for all $i =1,\ldots,n $.
\end{proposition}
\begin{proposition} \cite{ky}
If $b \in {\mathbb{H}}$, then  ${\rm{cdet}} _{{j}}\, {\rm {\bf
A}}_{.j} \left( {{\rm {\bf a}}_{.j} b} \right) = {\rm{cdet}}
_{{j}}\, {\rm {\bf A}}  b$ for all $j =1,\ldots,n$.
\end{proposition}
\begin{proposition} \cite{ky}
If for  ${\rm {\bf A}}\in {\rm M}\left( {n,{\mathbb{H}}}
\right)$\, there exists $t \in I_{n} $ such that $a_{tj} = b_{j} +
c_{j} $\, for all $j = 1,\ldots,n$, then
\[
\begin{array}{l}
   {\rm{rdet}}_{{i}}\, {\rm {\bf A}} = {\rm{rdet}}_{{i}}\, {\rm {\bf
A}}_{{t{\kern 1pt}.}} \left( {{\rm {\bf b}}} \right) +
{\rm{rdet}}_{{i}}\, {\rm {\bf A}}_{{t{\kern 1pt}.}} \left( {{\rm
{\bf c}}} \right), \\
  {\rm{cdet}} _{{i}}\, {\rm {\bf A}} = {\rm{cdet}} _{{i}}\, {\rm
{\bf A}}_{{t{\kern 1pt}.}} \left( {{\rm {\bf b}}} \right) +
{\rm{cdet}}_{{i}}\, {\rm {\bf A}}_{{t{\kern 1pt}.}} \left( {{\rm
{\bf c}}} \right),
\end{array}
\]
\noindent where ${\rm {\bf b}}=(b_{1},\ldots, b_{n})$, ${\rm {\bf
c}}=(c_{1},\ldots, c_{n})$ and for all $i =1,\ldots,n$.
\end{proposition}
\begin{proposition} \cite{ky}
If for ${\rm {\bf A}}\in {\rm M}\left( {n,{\mathbb{H}}} \right)$\,
 there exists $t \in J_{n} $ such that $a_{i\,t} = b_{i} + c_{i}$
for all $i = 1,\ldots,n$, then
\[
\begin{array}{l}
  {\rm{rdet}}_{{j}}\, {\rm {\bf A}} = {\rm{rdet}}_{{j}}\, {\rm {\bf
A}}_{{\,.\,{\kern 1pt}t}} \left( {{\rm {\bf b}}} \right) +
{\rm{rdet}}_{{j}}\, {\rm {\bf A}}_{{\,.\,{\kern 1pt} t}} \left(
{{\rm
{\bf c}}} \right),\\
  {\rm{cdet}} _{{j}}\, {\rm {\bf A}} = {\rm{cdet}} _{{j}}\, {\rm
{\bf A}}_{{\,.\,{\kern 1pt}t}} \left( {{\rm {\bf b}}} \right) +
{\rm{cdet}} _{{j}} {\rm {\bf A}}_{{\,.\,{\kern 1pt}t}} \left(
{{\rm {\bf c}}} \right),
\end{array}
\]
\noindent where ${\rm {\bf b}}=(b_{1},\ldots, b_{n})^T$, ${\rm
{\bf c}}=(c_{1},\ldots, c_{n})^T$ and for all $j =1,\ldots,n$.
\end{proposition}

The following lemmas enable us to expand ${\rm{rdet}}_{{i}}\, {\rm
{\bf A}}$ by cofactors
  along  the $i$th row and ${\rm{cdet}} _{j} {\rm {\bf A}}$
 along  the $j$th column respectively for all $i, j = 1,\ldots,n$.

\begin{lemma}\label{kyrc1}  \cite{ky}
Let $R_{i{\kern 1pt} j}$ be the right $ij$-th cofactor of ${\rm
{\bf A}}\in {\rm M}\left( {n, {\mathbb{H}}} \right)$, that is, $
{\rm{rdet}}_{{i}}\, {\rm {\bf A}} = {\sum\limits_{j = 1}^{n}
{{a_{i{\kern 1pt} j} \cdot R_{i{\kern 1pt} j} } }} $ for all $i =
1,\ldots,n$.  Then
\[
 R_{i{\kern 1pt} j} = {\left\{ {{\begin{array}{*{20}c}
  - {\rm{rdet}}_{{j}}\, {\rm {\bf A}}_{{.{\kern 1pt} j}}^{{i{\kern 1pt} i}} \left( {{\rm
{\bf a}}_{{.{\kern 1pt} {\kern 1pt} i}}}  \right),& {i \ne j},
\hfill \\
 {\rm{rdet}} _{{k}}\, {\rm {\bf A}}^{{i{\kern 1pt} i}},&{i = j},
\hfill \\
\end{array}} } \right.}
\]
\noindent where  ${\rm {\bf A}}_{.{\kern 1pt} j}^{i{\kern 1pt} i}
\left( {{\rm {\bf a}}_{.{\kern 1pt} {\kern 1pt} i}}  \right)$ is
obtained from ${\rm {\bf A}}$   by replacing the $j$th column with
the $i$th column, and then by deleting both the $i$th row and
column, $k = \min {\left\{ {I_{n}}  \right.} \setminus {\left.
{\{i\}} \right\}} $.
\end{lemma}
\begin{lemma}\label{kyrc5} \cite{ky}
Let $L_{i{\kern 1pt} j} $ be the left $ij$-th cofactor of
 ${\rm {\bf A}}\in {\rm M}\left( {n,{\mathbb{H}}} \right)$, that
 is,
$ {\rm{cdet}} _{{j}}\, {\rm {\bf A}} = {{\sum\limits_{i = 1}^{n}
{L_{i{\kern 1pt} j} \cdot a_{i{\kern 1pt} j}} }}$ for all $j
=1,\ldots,n$. Then
\[
L_{i{\kern 1pt} j} = {\left\{ {\begin{array}{*{20}c}
 -{\rm{cdet}} _{i}\, {\rm {\bf A}}_{i{\kern 1pt} .}^{j{\kern 1pt}j} \left( {{\rm {\bf a}}_{j{\kern 1pt}. } }\right),& {i \ne
j},\\
 {\rm{cdet}} _{k}\, {\rm {\bf A}}^{j\, j},& {i = j},
\\
\end{array} }\right.}
\]
\noindent where  ${\rm {\bf A}}_{i{\kern 1pt} .}^{jj} \left( {{\rm
{\bf a}}_{j{\kern 1pt} .} } \right)$ is obtained from ${\rm {\bf
A}}$
 by replacing the $i$th row with the $j$th row, and then by
deleting both the $j$th row and  column, $k = \min {\left\{
{J_{n}} \right.} \setminus {\left. {\{j\}} \right\}} $.
\end{lemma}
We recall some well-known definitions. The \textit{conjugate} of a
quaternion $a = a_{0} + a_{1} i + a_{2} j + a_{3} k \in {\rm
{\mathbb{H}}}$ is defined by $\overline {a} = a_{o}-a_{1}i- a_{2}j
- a_{3}k$. The \textit{Hermitian adjoint matrix} of ${\rm {\bf A}}
= \left( {a_{ij}} \right) \in {\rm {\mathbb{H}}}^{n\times m}$ is
called the matrix ${\rm {\bf A}}^{ *} = \left( {a_{ij}^{ *} }
\right)_{m\times n} $ if $a_{ij}^{ *} = \overline {a_{ji}}  $ for
all $i = 1,\ldots,n $ and $j = 1,\ldots,m$.
 The matrix ${\rm {\bf A}} = \left( {a_{ij}}  \right) \in {\rm
{\mathbb{H}}}^{n\times m}$ is \textit{Hermitian} if ${\rm {\bf
A}}^{ *}  = {\rm {\bf A}}$.

A following theorem has a key value in the theory of the column
and row determinants.

\begin{theorem} \cite{ky}\label{kyrc2}
If ${\rm {\bf A}} = \left( {a_{ij}}  \right) \in {\rm M}\left(
{n,{\rm {\mathbb{H}}}} \right)$ is Hermitian, then ${\rm{rdet}}
_{1} {\rm {\bf A}} = \cdots = {\rm{rdet}} _{n} {\rm {\bf A}} =
{\rm{cdet}} _{1} {\rm {\bf A}} = \cdots = {\rm{cdet}} _{n} {\rm
{\bf A}} \in {\rm {\mathbb{R}}}.$
\end{theorem}
 Taking into account Theorem \ref{kyrc2} we define the determinant of a
Hermitian matrix by putting $\det {\rm {\bf A}}: = {\rm{rdet}}
_{i} {\rm {\bf A}} = {\rm{cdet}} _{i} {\rm {\bf A}}$ for all $i
=1,\ldots,n $. This determinant of a Hermitian matrix coincides
with the Moore determinant. The properties of the determinant of a
Hermitian matrix are considered in \cite{ky} by means of the
column and row determinants. Among them we note the followings.

\begin{theorem} \cite{ky}
If the $i$th row of a Hermitian matrix ${\rm {\bf A}}\in {\rm
M}\left( {n,{\rm {\mathbb{H}}}} \right)$ is replaced with a left
linear combination of its other rows, i.e. ${\rm {\bf a}}_{i.} =
c_{1} {\rm {\bf a}}_{i_{1} .} + \cdots + c_{k}  {\rm {\bf
a}}_{i_{k} .}$, where $ c_{l} \in {\rm {\mathbb{H}}}$ for all $l
=1,\ldots,k$ and $\{i,i_{l}\}\subset I_{n} $, then
\[{\rm{cdet}}_{i} {\rm {\bf A}}_{i.} \left( {c_{1} \cdot {\rm {\bf
a}}_{i_{1} .} + \cdots + c_{k} \cdot {\rm {\bf a}}_{i_{k} .}}\,
\right) = {\rm{rdet}}_{i} {\rm {\bf A}}_{i.} \left( {c_{1} \cdot
{\rm {\bf a}}_{i_{1} .} + \cdots + c_{k} \cdot {\rm {\bf
a}}_{i_{k} .}}\, \right) = 0.\]
\end{theorem}

\begin{theorem} \cite{ky}
If the $j$th column of
 a Hermitian matrix
${\rm {\bf A}}\in {\rm M}\left( {n,{\rm {\mathbb{H}}}} \right)$ is
replaced with a right linear combination of its other columns,
i.e. ${\rm {\bf a}}_{.j} = {\rm {\bf a}}_{.j_{1}}   c_{1} + \cdots
+ {\rm {\bf a}}_{.j_{k}} c_{k} $, where $c_{l} \in {\rm
{\mathbb{H}}}$ for all $l =1,\ldots,k$ and $\{j,j_{l}\}\subset
J_{n}$,  then
\[{\rm{cdet}}_{i} {\rm {\bf A}}_{.\,i} \left( {{\rm {\bf
a}}_{.\,i_{1}} \cdot c_{1} + \cdots + {\rm {\bf a}}_{.\,i_{k}}
\cdot c_{k}} \right) = {\rm{rdet}}_{i} {\rm {\bf A}}_{.\,i} \left(
{{\rm {\bf a}}_{.\,i_{1}} \cdot c_{1} + \cdots + {\rm {\bf
a}}_{.\,i_{k}} \cdot c_{k} } \right) = 0.\]
\end{theorem}
The following theorem about determinantal representation of an
inverse matrix of Hermitian follows immediately from these
properties.

\begin{theorem}
 \cite{ky}\label{kyrc3}
 If a Hermitian matrix ${\rm
{\bf A}} \in {\rm M}\left( {n,{\rm {\mathbb{H}}}} \right)$ is such
that $\det {\rm {\bf A}} \ne 0$, then there exist a unique right
inverse  matrix $(R{\rm {\bf A}})^{ - 1}$ and a unique left
inverse matrix $(L{\rm {\bf A}})^{ - 1}$, and $\left( {R{\rm {\bf
A}}} \right)^{ - 1} = \left( {L{\rm {\bf A}}} \right)^{ - 1} =
:{\rm {\bf A}}^{ - 1}$. They possess the following determinantal
representations:
\begin{equation}
\label{kyr1}
  \left( {R{\rm {\bf A}}} \right)^{ - 1} = {\frac{{1}}{{\det {\rm
{\bf A}}}}}
\begin{pmatrix}
  R_{11} & R_{21} & \cdots & R_{n1} \\
  R_{12} & R_{22} & \cdots & R_{n2} \\
  \cdots & \cdots & \cdots & \cdots \\
  R_{1n} & R_{2n} & \cdots & R_{nn}
\end{pmatrix},
\end{equation}
  \begin{equation}
\label{kyr2} \left( {L{\rm {\bf A}}} \right)^{ - 1} =
{\frac{{1}}{{\det {\rm {\bf A}}}}}
\begin{pmatrix}
  L_{11} & L_{21} & \cdots  & L_{n1} \\
  L_{12} & L_{22} & \cdots  & L_{n2} \\
  \cdots  & \cdots  & \cdots  & \cdots  \\
  L_{1n} & L_{2n} & \cdots  & L_{nn}
\end{pmatrix}.
\end{equation}
Here $R_{ij}$,  $L_{ij}$ are right and left $ij$-th cofactors of
${\rm {\bf
 A}}$ respectively for all $i,j =1,\ldots,n$.
\end{theorem}
 To obtain
determinantal representation of an arbitrary inverse matrix ${\rm
{\bf A}}^{ - 1}$, we consider the right ${\rm {\bf A}}{\rm {\bf
A}}^{ *} $ and left ${\rm {\bf A}}^{ *} {\rm {\bf A}}$
corresponding Hermitian matrix.

\begin{theorem} \cite{ky}
If an arbitrary column  of ${\rm {\bf A}}\in {\rm
{\mathbb{H}}}^{m\times n} $ is a right  linear combination of its
other columns, or an arbitrary row of ${\rm {\bf A}}^{ * }$ is a
left linear combination of its others, then $\det {\rm {\bf A}}^{
* }{\rm {\bf A}} = 0.$
\end{theorem}
\noindent Since the principal submatrices of a Hermitian matrix
are Hermitian, the principal minor may be defined as the
determinant of its principal submatrix by analogy to the
commutative case. We introduce \textit{the rank by principal
minors} that is the maximal order of a nonzero principal minor of
a Hermitian matrix. The following theorem determines a
relationship between it and the rank of a matrix defining as
ceiling amount of right-linearly independent columns  or
left-linearly independent rows which form basis.

\begin{theorem} \cite{ky}
A rank by principal minors of ${\rm {\bf A}}^{ *} {\rm {\bf A}}$
is equal to its rank and a rank of ${\rm {\bf A}} \in {\rm
{\mathbb{H}}}^{m\times n}$.
\end{theorem}
\begin{theorem} \cite{ky}
If ${\rm {\bf A}} \in {\rm {\mathbb{H}}}^{m\times n}$, then an
arbitrary column of  ${\rm {\bf A}}$ is a right linear combination
of its basis columns or an arbitrary row of  ${\rm {\bf A}}$ is a
left linear combination of its basis rows.
\end{theorem}
\noindent The criterion of singularity of a Hermitian matrix is
obtained.

\begin{theorem} \cite{ky}
The right-linearly independence of columns of ${\rm {\bf A}} \in
{\rm {\mathbb{H}}}^{m\times n}$ or the left-linearly independence
of rows of ${\rm {\bf A}}^{ *} $ is the necessary and sufficient
condition for $\det {\rm {\bf A}}^{ *} {\rm {\bf A}} \neq 0.$
\end{theorem}
\begin{theorem}\cite{ky}
If ${\rm {\bf A}} \in {\rm M}\left( {n,{\rm {\mathbb{H}}}}
\right)$, then $\det {\rm {\bf A}}{\rm {\bf A}}^{ *} = \det {\rm
{\bf A}}^{ *} {\rm {\bf A}}$.
\end{theorem}
\noindent A concept of the double determinant is introduced by
this theorem. This concept was initially introduced by L. Chen in
\cite{ch}.
\begin{definition}
The determinant of the corresponding Hermitian matrix of ${\rm
{\bf A}} \in {\rm M}\left( {n,{\rm {\mathbb{H}}}} \right)$
 is called its double determinant, i.e.
${\rm{ddet}}{ \rm{\bf A}}: = \det \left( {{\rm {\bf A}}^{ *} {\rm
{\bf A}}} \right) = \det \left( {{\rm {\bf A}}{\rm {\bf A}}^{ *} }
\right)$.
\end{definition}
\noindent The relationship between the double determinant and the
noncommutative determinants of E. Moore, E. Study and J. Diedonne
is obtained, ${\rm ddet} {\rm {\bf A}} = {\rm Mdet} \left( {{\rm
{\bf A}}^{ *} {\rm {\bf A}}} \right) = {\rm Sdet} {\rm {\bf A}} =
{\rm Ddet} ^{2}{\rm {\bf A}}$. But unlike those, the double
determinant can be expanded  along an arbitrary row or column  by
means of the column and row determinants.

\begin{definition}
Suppose ${\rm {\bf A}} \in {\rm M}\left( {n,{\rm {\mathbb{H}}}}
\right)$. We have a column expansion of ${\rm{ddet}} {\rm {\bf
A}}$ along the $j$th column, ${\rm{ddet}} {\rm {\bf A}} =
{\rm{cdet}} _{j} \left( {{\rm {\bf A}}^{ *} {\rm {\bf A}}} \right)
= {\sum\limits_{i} {{\mathbb{L}} _{ij} \cdot a_{ij}} }$, and a row
expansion of it along the $i$th row,
 ${\rm{ddet}} {\rm {\bf A}} = {\rm{rdet}}_{i} \left(
{{\rm {\bf A}}{\rm {\bf A}}^{ *} } \right) = {\sum\limits_{j}
{a_{ij} \cdot} } {\mathbb{R}} _{ i{\kern 1pt}j} $ for all $i,j
=1,\ldots,n $. Then by definition of the left double $ij$th
cofactor we put ${\mathbb{L}} _{ij} $ and by definition of the
right double $ij$th cofactor we put ${\mathbb{R}} _{ i{\kern1pt}j}
$.
\end{definition}

\begin{theorem} \cite{ky}\label{kyrc4}
The necessary and sufficient condition of invertibility of  ${\rm
{\bf A}} \in {\rm M}(n,{\rm {\mathbb{H}}})$ is ${\rm{ddet}} {\rm
{\bf A}} \ne 0$. Then there exists $ {\rm {\bf A}}^{ - 1} = \left(
{L{\rm {\bf A}}} \right)^{ - 1} = \left( {R{\rm {\bf A}}}
\right)^{ - 1}$, where
\begin{equation}
\label{kyr3} \left( {L{\rm {\bf A}}} \right)^{ - 1} =\left( {{\rm
{\bf A}}^{ *}{\rm {\bf A}} } \right)^{ - 1}{\rm {\bf A}}^{ *}
={\frac{{1}}{{{\rm{ddet}}{ \rm{\bf A}} }}}
\begin{pmatrix}
  {\mathbb{L}} _{11} & {\mathbb{L}} _{21}& \ldots & {\mathbb{L}} _{n1} \\
  {\mathbb{L}} _{12} & {\mathbb{L}} _{22} & \ldots & {\mathbb{L}} _{n2} \\
  \ldots & \ldots & \ldots & \ldots \\
 {\mathbb{L}} _{1n} & {\mathbb{L}} _{2n} & \ldots & {\mathbb{L}} _{nn}
\end{pmatrix},
\end{equation}
\begin{equation}\label{kyr4} \left( {R{\rm {\bf A}}} \right)^{ - 1} = {\rm {\bf
A}}^{ *} \left( {{\rm {\bf A}}{\rm {\bf A}}^{ *} } \right)^{ - 1}
= {\frac{{1}}{{{\rm{ddet}}{ \rm{\bf A}} }}}
\begin{pmatrix}
 {\mathbb{R}} _{\,{\kern 1pt} 11} & {\mathbb{R}} _{\,{\kern 1pt} 21} &\ldots & {\mathbb{R}} _{\,{\kern 1pt} n1} \\
 {\mathbb{R}} _{\,{\kern 1pt} 12} & {\mathbb{R}} _{\,{\kern 1pt} 22} &\ldots & {\mathbb{R}} _{\,{\kern 1pt} n2}  \\
 \ldots  & \ldots & \ldots & \ldots \\
 {\mathbb{R}} _{\,{\kern 1pt} 1n} & {\mathbb{R}} _{\,{\kern 1pt} 2n} &\ldots & {\mathbb{R}} _{\,{\kern 1pt} nn}
\end{pmatrix},
\end{equation}
\noindent and ${\mathbb{L}} _{ij} = {\rm{cdet}} _{j} ({\rm {\bf
A}}^{\ast}{\rm {\bf A}})_{.j} \left( {{\rm {\bf a}}_{.{\kern 1pt}
i}^{ *} } \right)$, ${\mathbb{R}} _{\,{\kern 1pt} ij} =
{\rm{rdet}}_{i} ({\rm {\bf A}}{\rm {\bf A}}^{\ast})_{i.} \left(
{{\rm {\bf a}}_{j.}^{ *} }  \right)$  for all $i,j =1,\ldots,n$.
\end{theorem}
\noindent This theorem introduces the determinantal
representations of an inverse matrix by the left (\ref{kyr1}) and
right (\ref{kyr2}) double cofactors. The inverse  matrix ${\rm
{\bf A}}^{ - 1}$ of ${\rm {\bf A}} \in {\rm M}(n,{\mathbb H})$ on
the assumption of ${\rm{ddet}}{ \rm{\bf A}} \ne 0$ is represented
by the
 analog of the classical adjoint matrix. If we denote this analog of the
adjoint matrix over ${\mathbb H}$ by ${\rm Adj}{\left[ {{\left[
{{\rm {\bf A}}} \right]}} \right]}$, where ${\rm Adj}{\left[
{{\left[ {{\rm {\bf A}}} \right]}} \right]}= ({\mathbb{L}}
_{ij})_{n\times n}$ or ${\rm Adj}{\left[ {{\left[ {{\rm {\bf A}}}
\right]}} \right]}= ({\mathbb{R}} _{ij})_{n\times n}$, then the
following formula is valid over
 ${\mathbb H}$:
\[
{\rm {\bf A}}^{ - 1} = {\frac{{{\rm Adj}{\left[ {{\left[ {{\rm
{\bf A}}} \right]}} \right]}}}{{{\rm{ddet}} {\rm {\bf A}}}}}.
\]
Using the determinantal representations of an inverse matrix by
the left (\ref{kyr1}) and right (\ref{kyr2})  analogs of a
classical adjoint matrix we obtain the Cramer rule for  right and
left systems of  linear equations respectively.

\begin{theorem} Let
$
 {\rm {\bf A}} \cdot {\rm {\bf x}} = {\rm {\bf y}}
$ be a right system of  linear equations  with a matrix of
coefficients ${\rm {\bf A}}\in {\rm M}(n,{\mathbb H})$, a column
of constants ${\rm {\bf y}} = \left( {y_{1} ,\ldots ,y_{n} }
\right)^{T}\in {\mathbb H}^{n\times 1}$ and a column of unknowns
${\rm {\bf x}} = \left( {x_{1} ,\ldots ,x_{n}} \right)^{T}$. If
${\rm{ddet}}{ \rm{\bf A}} \ne 0$, then we have for all$j =
\overline {1,n}$
\[
 x_{j}^{} = {\frac{{{\rm{cdet}} _{j} ({\rm {\bf A}}^{
*} {\rm {\bf A}})_{.j} \left( {{\rm {\bf f}}}
\right)}}{{{\rm{ddet}} {\rm {\bf A}}}}},
\]
where $ {\rm {\bf f}} = {\rm {\bf A}}^{ *} {\rm {\bf y}}.$
\end{theorem}
\begin{theorem} Let
$
 {\rm {\bf x}} \cdot {\rm {\bf A}} = {\rm {\bf y}}
$ be a left  system of  linear equations
 with a
matrix of coefficients ${\rm {\bf A}}\in {\rm M}(n,{\mathbb H})$,
a row of constants ${\rm {\bf y}} = \left( {y_{1} ,\ldots ,y_{n} }
\right)\in {\mathbb H}^{1\times n}$ and a row of unknowns ${\rm
{\bf x}} = \left( {x_{1} ,\ldots ,x_{n}}  \right)$. If
${\rm{ddet}}{ \rm{\bf A}} \ne 0$, then we obtain for all $i =
\overline {1,n}$
\[
  x_{i}^{} = {\frac{{{\rm{rdet}}_{i} ({\rm {\bf A}}{\rm {\bf A}}^{ *})_{i.} \left( {{\rm
{\bf z}}} \right)}}{{{\rm{ddet}}{ \rm{\bf A}}}}},
\]
where $ {\rm {\bf z}} = {\rm {\bf y}}{\rm {\bf A}}^{ *}$.
\end{theorem}

\section{Cramer's rule for some matrix equations.}
We denote ${\rm {\bf A}}^{ \ast}{\rm {\bf B}}=:\hat{{\rm {\bf
B}}}= (\hat{b}_{ij})$, ${\rm {\bf B}}{\rm {\bf A}}^{
\ast}=:\check{{\rm {\bf B}}}= (\check{b}_{ij})$.
\begin{theorem} Suppose
\begin{equation}\label{kyr5}
 {\rm {\bf A}}{\rm {\bf X}} = {\rm {\bf B}}
\end{equation}
\noindent is a right matrix equation, where ${\left\{ {{\rm {\bf
A}},{\rm {\bf B}}} \right\}} \in {\rm M}(n,{\rm {\mathbb{H}}} )$
are given, ${\rm {\bf X}} \in {\rm M}(n,{\rm {\mathbb{H}}} )$ is
unknown. If ${\rm ddet} {\rm {\bf A}} \ne 0$, then (\ref{kyr5})
has a unique solution, and the solution is
\begin{equation}
\label{kyr6} x_{i\,j} = {\frac{{{\rm cdet} _{i} ({\rm {\bf A}}^{
*} {\rm {\bf A}})_{.\,i} \left( {\hat{{\rm {\bf b}}}_{.j}}
\right)}}{{ {\rm ddet} {\rm {\bf A}}}}}
\end{equation}
\noindent where $\hat{{\rm {\bf b}}}_{.j}$ is the $j$th column of
$\hat{{\rm {\bf B}}}$ for all $i,j = 1,...,n $.
\end{theorem}
{\textit{Proof.}} By Theorem \ref{kyrc4} the matrix ${\rm {\bf
A}}$ is invertible. There exists the unique inverse matrix ${\rm
{\bf A}}^{ - 1}$. From this it follows that the solution of
(\ref{kyr5}) exists and is unique, ${\rm {\bf X}} = {\rm {\bf
A}}^{ - 1}{\rm {\bf B}}$. If we represent ${\rm {\bf A}}^{ -
1}=({\rm {\bf A}}^{ \ast}{\rm {\bf A}})^{ - 1}{\rm {\bf A}}^{
\ast}$ as a left inverse by (\ref{kyr3}),  and use the
determinantal representation of $({\rm {\bf A}}^{ \ast}{\rm {\bf
A}})^{ - 1}$ by (\ref{kyr2}), then for all $i,j = 1,...,n $ we
obtain
\[
x_{ij}  = {\frac{{1}}{{{\rm ddet} {\rm {\bf A}}}}}{\sum\limits_{k
= 1}^{n} {{L}_{ki} \hat{b}_{kj}} },
\]
\noindent where ${L}_{ij}$ is a left $ij$th
 cofactor of $({\rm {\bf A}}^{ \ast}{\rm {\bf A}})$ for all $i,j = 1,...,n
$. From this by  Lemma \ref{kyrc5}  and denoting the $j$-th column
of $\hat{{\rm {\bf B}}}$ by $\hat{{\rm {\bf b}}}_{.j} $, it
follows (\ref{kyr6}).$\blacksquare$

\begin{theorem} Suppose
\begin{equation}\label{kyr7}
 {\rm {\bf X}}{\rm {\bf A}} = {\rm {\bf B}}
\end{equation}
\noindent is a left matrix equation, where ${\left\{ {{\rm {\bf
A}},{\rm {\bf B}}} \right\}} \in {\rm M}(n,{\rm {\mathbb{H}}} )$
are given, ${\rm {\bf X}} \in {\rm M}(n,{\rm {\mathbb{H}}} )$ is
unknown. If ${\rm ddet} {\rm {\bf A}} \ne 0$, then (\ref{kyr7})
has a unique solution, and the solution is
\begin{equation}
\label{kyr8} x_{i\,j} = {\frac{{{\rm rdet} _{j} ({\rm {\bf A}}{\rm
{\bf A}}^{ *} )_{j.\,} \left( {\check{{\rm {\bf b}}}_{i\,.}}
\right)}}{{{\rm ddet} {\rm {\bf A}}}}}
\end{equation}
\noindent where $\check{{\rm {\bf b}}}_{i.}$ is the $i$th column
of $\check{{\rm {\bf B}}}$ for all $i,j = 1,...,n $.
\end{theorem}
{\textit{Proof.}} By Theorem \ref{kyrc4} the matrix ${\rm {\bf
A}}$ is invertible. There exists the unique inverse matrix ${\rm
{\bf A}}^{ - 1}$. From this it follows that the solution of
(\ref{kyr7}) exists and is unique, ${\rm {\bf X}} = {\rm {\bf
B}}{\rm {\bf A}}^{ - 1}$. If we represent $\left( {{\rm {\bf A}}}
\right)^{ - 1} = {\rm {\bf A}}^{ *} \left( {{\rm {\bf A}}{\rm {\bf
A}}^{ *} } \right)^{ - 1}$ as a right inverse by (\ref{kyr4}) and
use the determinantal representation of $({\rm {\bf A}}{\rm {\bf
A}}^{ \ast})^{ - 1}$ by (\ref{kyr1}), then for all $i,j =1,...,n $
we have
\[
x_{ij} = {\frac{{1}}{{{\rm ddet} {\rm {\bf A}}}}}{\sum\limits_{k =
1}^{n} {\check{b}_{i\,k}} } { R}_{jk}.
\]

\noindent where ${R}_{i{\kern 1pt} j}$ is a right $ij$th cofactor
of $({\rm {\bf A}}{\rm {\bf A}}^{ \ast})$ for all $i,j = 1,...,n
$. From this by means of Lemma \ref{kyrc1} and denoting the $i$th
row of $\check{{\rm {\bf B}}}$ by $\check{{\rm {\bf b}}}_{i\,.} $,
it follows (\ref{kyr8}). $\blacksquare$

We denote ${\rm {\bf A}}^{ \ast}{\rm {\bf C}}{\rm {\bf B}}^{
\ast}=:\tilde{{\rm {\bf C}}}=(\tilde{c}_{ij})$.
\begin{theorem} Suppose
\begin{equation}\label{kyr9}
{\rm {\bf A}}{\rm {\bf X}}{\rm {\bf B}} = {\rm {\bf C}}
\end{equation}
\noindent is a two-sided matrix equation, where ${\left\{ {{\rm
{\bf A}},{\rm {\bf B}},{\rm {\bf C}}} \right\}} \in {\rm M}(n,{\rm
{\mathbb{H}}} )$ are given, ${\rm {\bf X}} \in {\rm M}(n,{\rm
{\mathbb{H}}} )$ is unknown. If ${\rm ddet} {\rm {\bf A}} \ne 0$
and ${\rm ddet} {\rm {\bf B}} \ne 0$ , then (\ref{kyr9}) has a
unique solution, and the solution is
\begin{equation}
\label{kyr10} x_{i\,j} = {\frac{{{\rm rdet} _{j} ({\rm {\bf
B}}{\rm {\bf B}}^{ *} )_{j.\,} \left( {{\rm {\bf c}}_{i\,.}^{{\rm
{\bf A}}}} \right)}}{{{\rm ddet}  {\rm {\bf A}}\cdot {\rm ddet}
{\rm {\bf B}}}}},
\end{equation}

\noindent or
\begin{equation}
\label{kyr11} x_{i\,j} = {\frac{{{\rm cdet} _{i} ({\rm {\bf A}}^{
*} {\rm {\bf A}})_{.\,i\,} \left( {{\rm {\bf c}}_{.j}^{{\rm {\bf
B}}}} \right)}}{{{\rm ddet} {\rm {\bf A}}\cdot {\rm ddet}  {\rm
{\bf B}}}}},
\end{equation}

\noindent where ${\rm {\bf c}}_{i\,.}^{{\rm {\bf A}}} : = \left(
{{\rm cdet} _{i} ({\rm {\bf A}}^{ *} {\rm {\bf A}})_{.\,i\,}
\left( {\tilde{{\rm {\bf c}}}_{.1}} \right),\ldots ,{\rm cdet}
_{i} ({\rm {\bf A}}^{ *} {\rm {\bf A}})_{.\,i\,} \left(
{\tilde{{\rm {\bf c}}}_{.n}} \right)} \right)$ is the row vector
and  ${\rm {\bf c}}_{.j}^{{\rm {\bf B}}} : = \left( {{\rm rdet}
_{j} ({\rm {\bf B}}{\rm {\bf B}}^{ *} )_{j.\,} \left( {\tilde{{\rm
{\bf c}}}_{1\,.}^{}} \right),\ldots ,{\rm rdet} _{j} ({\rm {\bf
B}}{\rm {\bf B}}^{ *} )_{j.\,} \left( {\tilde{{\rm {\bf
c}}}_{n\,.}^{}} \right)} \right)^{T}$ is the column vector and
$\tilde{{\rm {\bf c}}}_{i\,.}$, $\tilde{{\rm {\bf c}}}_{.\,j}$ are
the  ith row vector and the jth column vector of ${\rm {\bf
\widetilde{C}}}$, respectively, for all $i,j = 1,...,n $.
\end{theorem}
{\textit{Proof.}} By Theorem \ref{kyrc4} the matrices ${\rm {\bf
A}}$ and ${\rm {\bf B}}$ are invertible. There exist the unique
inverse matrices ${\rm {\bf A}}^{ - 1}$ and ${\rm {\bf B}}^{ -
1}$. From this it follows that the solution of (\ref{kyr9}) exists
and is unique, ${\rm {\bf X}} = {\rm {\bf A}}^{ - 1}{\rm {\bf
C}}{\rm {\bf B}}^{ - 1}$. If we represent ${\rm {\bf A}}^{ -
1}=({\rm {\bf A}}^{ \ast}{\rm {\bf A}})^{ - 1}{\rm {\bf A}}^{
\ast}$ as a left inverse  and $\left( {{\rm {\bf B}}} \right)^{ -
1} = {\rm {\bf B}}^{ *} \left( {{\rm {\bf B}}{\rm {\bf B}}^{ *} }
\right)^{ - 1}$ as a right inverse, then for all $i,j = 1,...,n $
we have
\[\begin{array}{c}
{\rm {\bf X}} = ({\rm {\bf A}}^{ \ast}{\rm {\bf A}})^{ - 1}{\rm
{\bf A}}^{ \ast}{\rm {\bf C}}{\rm {\bf B}}^{ *} \left( {{\rm {\bf
B}}{\rm {\bf B}}^{ *} } \right)^{ - 1}=\\
 = \begin{pmatrix}
   x_{11} & x_{12} & \ldots & x_{1n} \\
   x_{21} & x_{22} & \ldots & x_{2n} \\
   \ldots & \ldots & \ldots & \ldots \\
   x_{n1} & x_{n2} & \ldots & x_{nn} \
 \end{pmatrix}
 = {\frac{{1}}{{{\rm ddet} {\rm {\bf A}}}}}\begin{pmatrix}
  {L} _{11}^{{\rm {\bf A}}} & {L} _{21}^{{\rm {\bf A}}}&
   \ldots & {L} _{n1}^{{\rm {\bf A}}} \\
  {L} _{12}^{{\rm {\bf A}}} & {L} _{22}^{{\rm {\bf A}}} &
  \ldots & {L} _{n2}^{{\rm {\bf A}}} \\
  \ldots & \ldots & \ldots & \ldots \\
 {L} _{1n}^{{\rm {\bf A}}} & {L} _{2n}^{{\rm {\bf A}}}
  & \ldots & {L} _{nn}^{{\rm {\bf A}}}
\end{pmatrix}\times\\
  \times\begin{pmatrix}
    \tilde{c}_{11} & \tilde{c}_{12} & \ldots & \tilde{c}_{1n} \\
    \tilde{c}_{21} & \tilde{c}_{22} & \ldots & \tilde{c}_{2n} \\
    \ldots & \ldots & \ldots & \ldots \\
    \tilde{c}_{n1} & \tilde{c}_{n2} & \ldots & \tilde{c}_{nn} \
  \end{pmatrix}

{\frac{{1}}{{{\rm ddet} {\rm {\bf A}}}}}
\begin{pmatrix}
 {R} _{\, 11}^{{\rm {\bf B}}} & {R} _{\, 21}^{{\rm {\bf B}}}
 &\ldots & {R} _{\, n1}^{{\rm {\bf B}}} \\
 {R} _{\, 12}^{{\rm {\bf B}}} & {R} _{\, 22}^{{\rm {\bf B}}} &\ldots &
 {R} _{\, n2}^{{\rm {\bf B}}} \\
 \ldots  & \ldots & \ldots & \ldots \\
 {R} _{\, 1n}^{{\rm {\bf B}}} & {R} _{\, 2n}^{{\rm {\bf B}}} &
 \ldots & {R} _{\, nn}^{{\rm {\bf B}}}
\end{pmatrix},
\end{array}
\]
\noindent where ${L}_{ij}^{{\rm {\bf A}}} $ is a left $ij$th
 cofactor of $({\rm {\bf A}}^{ \ast}{\rm {\bf A}})$ and ${R}_{i{\kern 1pt} j}^{{\rm {\bf
B}}} $ is a right  $ij$th cofactor of $\left( {{\rm {\bf B}}{\rm
{\bf B}}^{ *} } \right)$ for all $i,j = 1,...,n $. This implies

\begin{equation}
\label{kyr12} x_{ij} = {\frac{{{\sum\limits_{m = 1}^{n} {\left(
{{\sum\limits_{k = 1}^{n} {{L}_{ki}^{{\rm {\bf A}}}
\tilde{c}_{\,km}} } } \right)}} {R}_{jm}^{{\rm {\bf B}}}} }{{{\rm
ddet} {\rm {\bf A}}\cdot{\rm ddet} {\rm {\bf B}}}}},
\end{equation}
\noindent for all $i,j = \overline {1,n} $. From this by Lemma
\ref{kyrc5}, we obtain
\[{\sum\limits_{k = 1}^{n}
{{L}_{ki}^{{\rm {\bf A}}} \tilde{c}_{k\,m}} }  = {\rm cdet} _{i}
({\rm {\bf A}}^{ *} {\rm {\bf A}})_{.\,i\,} \left( {\tilde{{\rm
{\bf c}}}_{.\,m}} \right),\] \noindent where $\tilde{{\rm {\bf
c}}}_{.\,m} $ is the $m$th column-vector of $\tilde{{\rm {\bf
C}}}$ for all $m = 1,...,n $. Denote by ${\rm {\bf
c}}_{i\,.}^{{\rm {\bf A}}} : = \left( { {\rm cdet} _{i} ({\rm {\bf
A}}^{ *} {\rm {\bf A}})_{.\,i\,} \left( {\tilde{{\rm {\bf
c}}}_{.1}} \right),\ldots , {\rm cdet} _{i} ({\rm {\bf A}}^{ *}
{\rm {\bf A}})_{.\,i} \left( {\tilde{{\rm {\bf c}}}_{.n}} \right)}
\right)$  the row-vector for all $i = 1,...,n $. Reducing the sum
${{\sum\limits_{m = 1}^{n} {\left( {{\sum\limits_{k = 1}^{n}
{{L}_{ki}^{{\rm {\bf A}}} \tilde{c}_{\,km}} } } \right)}}
{R}_{jm}^{{\rm {\bf B}}}} $ by Lemma \ref{kyrc1}, we obtain an
analog of Cramer's rule for (\ref{kyr9}) by (\ref{kyr10}).

Having changed the order of summation in (\ref{kyr12}), we have
\[
x_{ij} = {\frac{{{\sum\limits_{k = 1}^{n} {{ L}_{ki}^{{\rm {\bf
A}}}} }\left( {{\sum\limits_{m = 1}^{n} {\tilde{c}_{\,km}} }
{R}_{jm}^{{\rm {\bf B}}}} \right)}}{{{\rm ddet} {\rm {\bf
A}}\cdot{\rm ddet}{\rm {\bf B}}}}}.\]

By Lemma \ref{kyrc1}, we obtain ${\sum\limits_{m = 1}^{n}
{c_{k\,m}} } {R}_{j\,m}^{{\rm {\bf B}}} = {\rm rdet} _{j} ({\rm
{\bf B}}{\rm {\bf B}}^{ *} )_{j.\,} \left( {\tilde{{\rm {\bf
c}}}_{k\,.}} \right)$, where $\tilde{{\rm {\bf c}}}_{k\,.} $ is a
$k$th row-vector of $\tilde{{\rm {\bf C}}}$ for all $k = 1,...,n
$. We denote by
\[{\rm {\bf c}}_{.\,j}^{{\rm {\bf B}}} : = \left(
{{\rm rdet} _{j} ({\rm {\bf B}}{\rm {\bf B}}^{ *} )_{j.\,} \left(
{\tilde{{\rm {\bf c}}}_{1\,.}} \right),\ldots ,{\rm rdet} _{j}
({\rm {\bf B}}{\rm {\bf B}}^{ *} )_{j.\,} \left( {\tilde{{\rm {\bf
c}}}_{n\,.}^{}} \right)} \right)^{T}\]
 the column-vector for all
$j = 1,...,n $. Reducing the sum ${\sum\limits_{k = 1}^{n} {{
L}_{ki}^{{\rm {\bf A}}}} } \left( {{\sum\limits_{m = 1}^{n}
{\tilde{c}_{\,km}} }{ R}_{jm}^{{\rm {\bf B}}}} \right)$ by Lemma
\ref{kyrc5}, we obtain Cramer's rule for (\ref{kyr9}) by
(\ref{kyr11}). $\blacksquare$

In solving the matrix equations by Cramer's rules (\ref{kyr6}),
(\ref{kyr8}), (\ref{kyr10}), (\ref{kyr11}) we do not use the
complex representation of quaternion matrices and work only in the
quaternion skew field.
\section{Example}
Let us consider the two-sided matrix equation
\begin{equation}\label{kyr13}
{\rm {\bf A}}{\rm {\bf X}}{\rm {\bf B}} = {\rm {\bf C}}
\end{equation}
\noindent where ${\rm {\bf A}}=\begin{pmatrix}
  i & -j & k \\
  k & -i & 1 \\
  2 & k & -j
\end{pmatrix}$, ${\rm {\bf B}}=\begin{pmatrix}
  -k & j & 2 \\
 i & k & i \\
 -j & 1 & i
\end{pmatrix}$ and ${\rm {\bf C}}=\begin{pmatrix}
  1 & i & j \\
 k & j & -2 \\
 i & 1 & j
\end{pmatrix}$.  Then we have ${\rm {\bf A}}^{\ast}=\begin{pmatrix}
  -i & -k & 2 \\
  j & i & -k \\
  -k & 1 & j
\end{pmatrix}$, ${\rm {\bf A}}^{*}{\rm {\bf A}}=\begin{pmatrix}
  6 & j+3k & -j-k \\
  -j-3k & 3 & i \\
  j+k & -i & 3
\end{pmatrix}$ and

${\rm {\bf B}}^{\ast}=\begin{pmatrix}
  k & -i & j \\
 -j & -k & 1 \\
 2 & -i & -i
\end{pmatrix}$,
${\rm {\bf B}}{\rm {\bf B}}^{*}=\begin{pmatrix}
  6 & -3i+j & -i+j \\
 3i-j & 3 & 1+2k \\
 i-j & 1-2k & 3
\end{pmatrix}$,

$\tilde{{\rm {\bf C}}}={\rm {\bf A}}^{*}{\rm {\bf C}}{\rm {\bf
B}}^{\ast}=\begin{pmatrix}
  2k & 1-i-k & 3+i+3k \\
 -2-4i & -2+i-k & i-k \\
 -4+2i & 1+2i+j & 1+4i+j
\end{pmatrix}$. It is easy to get, ${\rm ddet} {\rm {\bf A}}=\det{\rm {\bf A}}^{*}{\rm {\bf
A}}=8$ and ${\rm ddet} {\rm {\bf B}}=\det{\rm {\bf B}}{\rm {\bf
B}}^{*}=4$. Therefore (\ref{kyr13}) has a solution. We shall find
it by (\ref{kyr10}). At first we obtain the row-vectors ${\rm {\bf
c}}_{i\,.}^{{\rm {\bf A}}}$ for all $i=1,2,3$.

\[\begin{array}{c}
    {\rm cdet} _{1} ({\rm {\bf A}}^{ *} {\rm {\bf
A}})_{.\,1\,} \left( {\tilde{{\rm {\bf c}}}_{.1}} \right)={\rm
cdet} _{1}\begin{pmatrix}
  2k & j+3k & -j-k \\
  -2-4i & 3 & i \\
  -4+2i & -i & 3
\end{pmatrix}=3\cdot 3 (2k)-\\-i(-i)(3j+5k)
  +(-j-k)(-i)(-2-4i)-3(j+3k)(-2-4i)+\\+(j+3k)i(-4+2i)-3(-j-k)(-4+2i)=\\
  =24j+8k,
\end{array}
\]
\noindent and so forth. Continuing in the same way, we get
\[\begin{array}{c}
  {\rm{\bf c}}_{1\,.}^{{\rm {\bf A}}}=(24j+8k, -8-8i+4j+4k,
8+8i+4j+4k), \\
   {\rm{\bf c}}_{2\,.}^{{\rm {\bf A}}}=(-20-36i, -10-2i-12j-12k,
-2-2i+12j+4k),\\
  {\rm{\bf c}}_{3\,.}^{{\rm {\bf A}}}=(12+4i, 6+2i+12j-4k,
6+10i-4j+4k).
\end{array}
\]
Then by (\ref{kyr10}) we have
\[\begin{array}{c}
  x_{1\,1} = {\frac{{{\rm rdet} _{1} ({\rm {\bf B}}{\rm {\bf B}}^{
*} )_{1.\,} \left( {{\rm {\bf c}}_{1\,.}^{{\rm {\bf A}}}}
\right)}}{{{\rm ddet}  {\rm {\bf A}}\cdot {\rm ddet}  {\rm {\bf
B}}}}}= \\
  \\= {\frac{1}{32}}\cdot {\rm rdet} _{1}\begin{pmatrix}
 24j+8k & -8-8i+4j+4k & 8+8i+4j+4k \\
 3i-j & 3 & 1+2k \\
 i-j & 1-2k & 3
\end{pmatrix}=\\
={\frac{1}{30}}\cdot ((24j+8k)\cdot 3\cdot
3-(24j+8k)(1+2k)(1-2k)+\\
+(-8-8i+4j+4k)(1+2k)(i-j)-(-8-8i+4j+4k)(3i-j)3+\\
+(8+8i+4j+4k)(1-2k)(3i-j)-(8+8i+4j+4k)(i-j)3)=\\=
{\frac{1}{32}}\cdot (-32+32i),
\end{array}
\]
\noindent and so forth. Continuing in the same way, we obtain
\[\begin{array}{c}
  x_{11}={\frac{-32+32i}{32}},
  x_{12}={\frac{-88-72i+24j-8k}{32}},
  x_{13}={\frac{24+8i-40j+56k }{32}}, \\
  x_{21}={\frac{-16i+32j-48k}{32}},
  x_{22}={\frac{20-28i-116j-76k}{32}},
  x_{23}={\frac{-44+68i+20j+12k}{32}}, \\
  x_{31}={\frac{16+16j+32k}{32}},
  x_{32}={\frac{20+44i+52j-28k}{32}},
  x_{32}={\frac{-12-20i+12j-4k}{32}}.
\end{array}
\]

\end{document}